\documentclass[12pt]{article}
\usepackage{amssymb,graphicx,rawfonts,t1enc}
\input{prepictex}
\input{pictex}
\input{postpictex}
\newcommand{\D}{\textbf{\textit{d}}}
\newcommand{\ZERO}{\textbf{0}}
\newcommand{\ONE}{\textbf{1}}

\newcommand{\wt}{\widetilde}
\newcommand{\F}{\mathbb F}
\newcommand{\R}{\mathbb R}
\newcommand{\C}{\mathbb C}
\newcommand{\Q}{\mathbb Q}
\begin{document}
\hfill
\vskip 3.0truecm
\thispagestyle{empty}
\par
\noindent
\centerline{{\large The Beckman-Quarles theorem for mappings from ${\R}^2$ to ${\F}^2$,}}
\par
\centerline{{\large where $\F$ is a subfield of a commutative field extending $\R$}}
\vskip 3.0truecm
\centerline{{\large Apoloniusz Tyszka}}
\vskip 3.0truecm
\par
\noindent
{\bf Summary.} Let $\F$ be a subfield of a commutative field extending $\R$.
Let $\varphi_2:{\F}^2 \times {\F}^2 \to~\F$,
$\varphi_2((x_1,x_2),(y_1,y_2))=(x_1-y_1)^2+(x_2-y_2)^2$.
We say that $f:{\R}^2 \to {\F}^2$ preserves distance $d \geq 0$ if
for each $x,y \in {\R}^2$ $|x-y|=d$ implies
$\varphi_2(f(x),f(y))=d^2$.
We prove that each unit-distance
preserving mapping $f:{\R}^2 \to {\F}^2$ has
a form $I \circ (\rho,\rho)$, where $\rho: \R \to \F$ is
a field homomorphism and $I: {\F}^2 \to {\F}^2$ is
an affine mapping with orthogonal linear part.
\vskip 3.0truecm
\par
Let $\F$ be a subfield of a commutative field extending $\R$.
Let $\varphi_n: {\F}^n \times {\F}^n \to \F$,
$\varphi_n((x_1,...,x_n),(y_1,...,y_n))=(x_1-y_1)^2+...+(x_n-y_n)^2$.
We say that $f:{\R}^n \to {\F}^n$ preserves distance $d \geq 0$ if
for each $x,y \in {\R}^n$ $|x-y|=d$ implies
$\varphi_n(f(x),f(y))=d^2$.
In this paper we study unit-distance preserving mappings from
${\R}^2$ to ${\F}^2$. We need the following technical Propositions 1-4.
\newpage
\par
\noindent
{\bf Proposition 1}
(\cite{Blumenthal}, \cite{Borsuk}, \cite{Tyszka200?}, \cite{Tyszka200??}).
The points 
$c_1,c_2,c_3 \in {\F}^2$
are affinely dependent if and only if their Cayley-Menger determinant
$\Delta(c_1,c_2,c_3):=$
$$
\det \left[
\begin{array}{cccc}
 0  &  1                 &  1                 & 1                  \\
 1  &  0                 & \varphi_2(c_1,c_2) & \varphi_2(c_1,c_3) \\
 1  & \varphi_2(c_2,c_1) &  0                 & \varphi_2(c_2,c_3) \\
 1  & \varphi_2(c_3,c_1) & \varphi_2(c_3,c_2) & 0                  \\
\end{array}\;\right]
$$
\par
\noindent
equals $0$.
\vskip 0.2truecm
\par
\noindent
{\bf Proposition 2}
(\cite{Blumenthal},\cite{Borsuk},\cite{Tyszka200?},\cite{Tyszka200??}).
For any points $c_1,c_2,c_3,c_4 \in {\F}^2$ their
Cayley-Menger determinant $\Delta(c_1,c_2,c_3,c_4):=$
$$
\det \left[
\begin{array}{ccccc}
 0  &  1                 &  1                 & 1                  & 1                  \\
 1  &  0                 & \varphi_2(c_1,c_2) & \varphi_2(c_1,c_3) & \varphi_2(c_1,c_4) \\
 1  & \varphi_2(c_2,c_1) &  0                 & \varphi_2(c_2,c_3) & \varphi_2(c_2,c_4) \\
 1  & \varphi_2(c_3,c_1) & \varphi_2(c_3,c_2) & 0                  & \varphi_2(c_3,c_4) \\
 1  & \varphi_2(c_4,c_1) & \varphi_2(c_4,c_2) & \varphi_2(c_4,c_3) & 0                  \\
\end{array}\;\right]
$$
\par
\noindent
equals $0$.
\vskip 0.2truecm
\par
\noindent
{\bf Proposition 3} (\cite{Tyszka200??}).
If $a,b \in \F$, $a+b \neq 0$, $z,x,\wt{x} \in {\F}^2$
and $\varphi_2(z,x)=a^2$, $\varphi_2(x,\wt{x})=~b^2$,
$\varphi_2(z,\wt{x})=(a+b)^2$, then
$\overrightarrow{zx}=\frac{\textstyle a}{{\textstyle a+b}}\overrightarrow{z\wt{x}}$.
\vskip 0.2truecm
\par
\noindent
{\it Proof.} Since $\varphi_2(z, \wt{x})=(a+b)^2 \neq 0$ we conclude that
$z \neq \wt{x}$.
The Cayley-Menger determinant $\Delta(z,x,\wt{x})=$
$$
\det \left[
\begin{array}{cccc}
0 &   1     &  1   &   1     \\
1 &   0     & a^2  & (a+b)^2 \\
1 &  a^2    &  0   &   b^2   \\
1 & (a+b)^2 & b^2  &   0     \\
\end{array}
\right]
=0,
$$
\par
\noindent
so by Proposition 1 the points $z,x,\wt{x}$ are affinely
dependent. Therefore, there exists $c \in \F$ such that
$\overrightarrow{zx}=c \cdot \overrightarrow{z\wt{x}}$.
Hence $a^2=\varphi_2(z,x)=c^2 \cdot\varphi_2(z,\wt{x})=c^2 \cdot (a+b)^2$.
Thus $c=\frac{\textstyle a}{\textstyle a+b}$
or $c=-\frac{\textstyle a}{\textstyle a+b}$ .
If $c=-\frac{\textstyle a}{\textstyle a+b}$ then
$\overrightarrow{x\wt{x}}=\overrightarrow{z\wt{x}}-\overrightarrow{zx}
=\overrightarrow{z\wt{x}}+\frac{\textstyle a}{\textstyle a+b}\overrightarrow{z\wt{x}}=
\frac{\textstyle 2a+b}{\textstyle a+b}\overrightarrow{z\wt{x}}$.
Hence $b^2=\varphi_2(x,\wt{x})=
\left( \frac{\textstyle 2a+b}{\textstyle a+b} \right)^2 \cdot \varphi_2(z,\wt{x})=
\left( \frac{\textstyle 2a+b}{\textstyle a+b} \right)^2 \cdot (a+b)^2=(2a+b)^2$.
Therefore $0=(2a+b)^2-b^2=4a(a+b)$. Since $a+b \neq 0$ we conclude that
$a=0$, so $c=-\frac{\textstyle a}{\textstyle a+b}=
\frac{\textstyle a}{\textstyle a+b}$ and the proof is completed.
\vskip 0.2truecm
\par
\noindent
{\bf Proposition 4.} If $E,F,C,D \in {\F}^2$, $\varphi_2(E,F) \neq 0$,
$C \neq D$ and
$\varphi_2(E,C)=\varphi_2(F,C)=\varphi_2(E,D)=\varphi_2(F,D)$,
then $\overrightarrow{EC}=\overrightarrow{DF}$ and
$\overrightarrow{FC}=\overrightarrow{DE}$ (see the points
$E,F,C,D$ in Figure 1).
\vskip 0.2truecm
\par
\noindent
{\it Proof.} Let $E=(e_1,e_2)$, $F=(f_1,f_2)$, $C=(c_1,c_2)$,
$D=(d_1,d_2)$. The mapping
$$
{\F}^2 \ni X \to X-\frac{E+F}{2} \in {\F}^2
$$
\par
\noindent
preserves both the values of $\overrightarrow{XY}$ for
all $X,Y \in {\F}^2$ and the values of $\varphi_2(X,Y)$
for all $X,Y \in {\F}^2$. Therefore, without loss of generality
we may assume that $(f_1,f_2)=F=-E=(-e_1,-e_2)$. Then
$$
(e_1-c_1)^2+(e_2-c_2)^2=\varphi_2(E,C)=\varphi_2(F,C)=(-e_1-c_1)^2+(-e_2-c_2)^2
$$
and 
$$
(e_1-d_1)^2+(e_2-d_2)^2=\varphi_2(E,D)=\varphi_2(F,D)=(-e_1-d_1)^2+(-e_2-d_2)^2.
$$
Hence
$$
~~~~~~~~~~~~~~~~~~~~~~~~~~~~~~~~~~~~~~~~~~~
e_1 \cdot c_1+e_2 \cdot c_2=0
~~~~~~~~~~~~~~~~~~~~~~~~~~~~~~~~~~~~~~~~~~~(1)
$$
and
$$
~~~~~~~~~~~~~~~~~~~~~~~~~~~~~~~~~~~~~~~~~~~
e_1 \cdot d_1+e_2 \cdot d_2=0
~~~~~~~~~~~~~~~~~~~~~~~~~~~~~~~~~~~~~~~~~~~(2).
$$
It implies that
$${e_1}^2+{c_1}^2+{e_2}^2+{c_2}^2=(e_1-c_1)^2+(e_2-c_2)^2=\varphi_2(E,C)=$$
$$\varphi_2(E,D)=(e_1-d_1)^2+(e_2-d_2)^2={e_1}^2+{d_1}^2+{e_2}^2+{d_2}^2.$$
Therefore
$$
~~~~~~~~~~~~~~~~~~~~~~~~~~~~~~~~~~~~~~~~~~~
{c_1}^2+{c_2}^2={d_1}^2+{d_2}^2
~~~~~~~~~~~~~~~~~~~~~~~~~~~~~~~~~~~~~~~~~~~(3).
$$
Since $0 \neq \varphi_2(E,F)=4 \cdot ({e_1}^2+{e_2}^2)$ we conclude that $e_1 \neq 0$
or $e_2 \neq 0$. Assume that $e_1 \neq 0$. Then by (1) and (2)
$c_1=\frac{{\textstyle -e_2 \cdot c_2}}{{\textstyle e_1}}$ and
$d_1=\frac{{\textstyle -e_2 \cdot d_2}}{{\textstyle e_1}}$, so applying (3)
we get
$$
\frac{{e_2}^2 \cdot {c_2}^2}{{e_1}^2}+{c_2}^2={c_1}^2+{c_2}^2={d_1}^2+{d_2}^2=
\frac{{e_2}^2 \cdot {d_2}^2}{{e_1}^2}+{d_2}^2.$$
It implies that
$$\varphi_2(E,F) \cdot {c_2}^2=4 \cdot ({e_1}^2+{e_2}^2) \cdot {c_2}^2=4 \cdot ({e_1}^2+{e_2}^2) \cdot {d_2}^2=\varphi_2(E,F) \cdot {d_2}^2.$$
Since $\varphi_2(E,F) \neq 0$ we conclude that ${c_2}^2={d_2}^2$.
It implies that $c_2=d_2$ or $c_2=-d_2$. In the fist case
$$c_1=\frac{{-e_2} \cdot {c_2}}{e_1}=\frac{{-e_2} \cdot {d_2}}{e_1}=d_1,$$
so $C=D$. Since $C \neq D$ we conclude that $c_2=-d_2$. Hence
$$c_1=\frac{{-e_2} \cdot {c_2}}{e_1}=\frac{\textstyle{{e_2} \cdot {d_2}}}{\textstyle{e_1}}=-d_1,$$
so $C=-D$. Analogical reasoning in the case $e_2 \neq 0$
gives also $C=-D$. Since $F=-E$ it proves that
$\overrightarrow{EC}=C-E=F-D=\overrightarrow{DF}$ and
$\overrightarrow{FC}=C-F=E-D=\overrightarrow{DE}$.
\vskip 0.6truecm
\par
If $\rho:\R \to \F$ is a field homomorphism then
$(\rho,...,\rho): {\R}^n \to {\F}^n$ preserves all distances
$\sqrt{r}$ with rational $r \geq 0$.
Indeed, if $(x_1-y_1)^2+...+(x_n-y_n)^2=(\sqrt{r})^2$ then
\begin{eqnarray*}
\varphi_n((\rho,...,\rho)(x_1,...,x_n),(\rho,...,\rho)(y_1,...,y_n))&=&\\
\varphi_n((\rho(x_1),...,\rho(x_n)),(\rho(y_1),...,\rho(y_n)))&=&\\
(\rho(x_1)-\rho(y_1))^2+...+(\rho(x_n)-\rho(y_n))^2&=&\\
(\rho(x_1-y_1))^2+...+(\rho(x_n-y_n))^2&=&\\
\rho((x_1-y_1)^2)+...+\rho((x_n-y_n)^2)&=&\\
\rho((x_1-y_1)^2+...+(x_n-y_n)^2)&=&\rho((\sqrt{r})^2)=\rho(r)=r=(\sqrt{r})^2.
\end{eqnarray*}
\vskip 0.2truecm
\par
Let $A_{n}(\F)$ denote the set of all positive numbers $d$
such that any map $f:{\R}^n \to {\F}^n$ that preserves
unit distance preserves also distance $d$.
The classical Beckman-Quarles theorem states that each unit-distance
preserving mapping from ${\R}^n$ to ${\R}^n$ ($n \geq 2$) is an isometry,
see \cite{Beckman-Quarles}-\cite{Benz1994} and \cite{Everling}.
It means that for each $n \geq 2$ $A_n(\R)=(0,\infty)$.
Let $D_{n}(\F)$ denote the set of all positive numbers~$d$
with the following property:
\vskip 0.2truecm
\par
\noindent
if $x,y \in {\R}^n$ and $|x-y|=d$ then there exists a finite set $S_{xy}$
with $\left\{x,y \right\} \subseteq S_{xy} \subseteq {\R}^n$ such that any map
$f:S_{xy}\rightarrow {\F}^n$ that preserves unit distance
preserves also the distance between $x$ and $y$.
\vskip 0.2truecm
\par
Obviously, $\{1\} \subseteq D_n(\F) \subseteq A_n(\F)$.
It is known 
(\cite{Tyszka200??}, reference \cite{Tyszka200?} contains a weaker result)
that $\{d>0: d^2 \in \Q\} \subseteq D_2(\F)$. In particular, all positive
rational numbers belong to $D_2(\F)$. Therefore, each unit-distance
preserving mapping $f:{\R}^2 \to {\F}^2$ preserves all rational
distances.
\vskip 0.2truecm
\par
\noindent
{\bf Theorem 1.} If $f:{\R}^2 \to {\F}^2$
preserves unit distance, then $f$ is injective.
\vskip 0.2truecm
\par
\noindent
{\it Proof.} Since $D_2(\F)$ is a dense subset of $(0,\infty)$,
for any $x,y \in {\R}^2$, $x \neq y$ there exists $z \in {\R}^2$
such that $|z-x| \neq |z-y|$ and $|z-x|,|z-y| \in D_2(\F)$.
All distances in $D_2(\F)$ are preserved by $f$.
Suppose $f(x)=f(y)$. This would imply that
$|z-x|^2=\varphi_2(f(z),f(x))=\varphi_2(f(z),f(y))=|z-y|^2$,
which is a contradiction.
\vskip 0.2truecm
\par
\noindent
{\bf Theorem 2.} If $f:{\R}^2 \to {\F}^2$ preserves unit distance,
$x,y \in {\R}^2$ and $x \neq y$, then $\varphi_2(f(x),f(y)) \neq 0$.
\vskip 0.2truecm
\par
\noindent
{\it Proof.} By Theorem 1 $f(x) \neq f(y)$. Since
$D_2(\F)$ is unbounded from above, there exist
$0<d \in D_2(\F)$ and $z \in {\R}^2$ such that $|x-z|=|y-z|=d$.
All distances in $D_2(\F)$ are preserved by $f$, so
$d^2=|x-z|^2=\varphi_2(f(x),f(z))$
and
$d^2=|y-z|^2=\varphi_2(f(y),f(z))$.
Assume, on the contrary, that $\varphi_2(f(x),f(y))=0$.
The Cayley-Menger determinant
$$\Delta(f(x),f(y),f(z))=
\det \left[
\begin{array}{cccc}
 0   &   1    &   1    &  1     \\
\ONE &  \ZERO &  \ZERO & {\D}^2 \\
\ONE &  \ZERO &  \ZERO & {\D}^2 \\
 1   &   d^2  &  d^2   &  0     \\
\end{array}\right]=0,
$$
so by Proposition 1 the points $f(x)$, $f(y)$, $f(z)$
are affinely dependent.
Thus there exists $c \in \F$ such that
$\overrightarrow{f(x)f(z)}=c\overrightarrow{f(x)f(y)}$.
Hence $d^2=\varphi_2(f(x),f(z))=c^2 \cdot \varphi_2(f(x),f(y))=0$,
which is a contradiction.
\vskip 0.2truecm
\par
\noindent
{\bf Theorem 3.} If $t \in \Q$, $0<t<1$, $A,B \in {\R}^2$, $A \neq B$,
$C=tA+(1-t)B$ and $f:{\R}^2 \to {\F}^2$ preserves unit distance,
then $f(C)=tf(A)+(1-t)f(B)$.
\vskip 0.2truecm
\par
\noindent
{\it Proof.} We choose $r \in \Q$ such that $r>|AB|$ and
$r<\frac{|AB|}{|1-2t|}$ if $t \neq \frac{1}{2}$. There exists
$D \in {\R}^2$ such that $|AD|=(1-t)r$ and $|BD|=tr$. Let
$E=tA+(1-t)D$ and $F=(1-t)B+tD$, see Figure 1.
\vskip 0.2truecm
\par
\centerline{
\beginpicture
\setcoordinatesystem units <0.8mm, 0.8mm>
\setplotsymbol({.})
\plot -20 0 40 0 /
\plot -6.67 25.82 40 0 /
\plot -23.33 12.9 -20 0 /
\arrow <4mm> [0.1,0.3] from -23.33 12.9 to 0 0
\arrow <4mm> [0.1,0.3] from -6.67 25.82 to 0 0
\arrow <4mm> [0.1,0.3] from -30 38.73 to -23.33 12.9
\arrow <4mm> [0.1,0.3] from -30 38.73 to -6.67 25.82
\put {$C$} at 0 -3
\put {$A$} at -20 -3
\put {$B$} at 40 -3
\put {$D$} at -30 41.73
\put {$F$} at -4.67 28.82
\put {$E$} at -26.33 12.9
\endpicture}
\vskip 0.2truecm
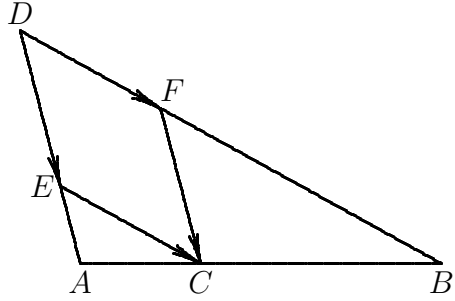
\centerline{Figure 1}
\vskip 0.2truecm
\par
\noindent
The segments $AE$, $ED$, $AD$, $BF$, $FD$, $BD$, $EC$
and $FC$ have rational lengths and
$|EC|=|FC|=|ED|=|FD|=t(1-t)r$.
Since $f$ preserves rational distances:
\\
\centerline{$\varphi_2(f(A),f(E))=|AE|^2=((1-t)^2r)^2,$}
\centerline{$\varphi_2(f(E),f(D))=|ED|^2=(t(1-t)r)^2,$}
\centerline{$\varphi_2(f(A),f(D))=|AD|^2=((1-t)r)^2.$}
Since $(1-t)^2r+t(1-t)r=(1-t)r$, by Proposition 3
$f(E)=tf(A)+(1-t)f(D)$. Analogously
$f(F)=(1-t)f(B)+tf(D)$. Since $C \neq D$, by Theorem 1
$f(C) \neq f(D)$.
Since $E \neq F$, by Theorem 2
$\varphi_2(f(E),f(F)) \neq 0$. Since $f$ preserves
rational distances:
\vskip 0.2truecm
\centerline{$\varphi_2(f(E),f(C))=|EC|^2=(t(1-t)r)^2,$}
\centerline{$\varphi_2(f(F),f(C))=|FC|^2=(t(1-t)r)^2,$}
\centerline{$\varphi_2(f(E),f(D))=|ED|^2=(t(1-t)r)^2,$}
\centerline{$\varphi_2(f(F),f(D))=|FD|^2=(t(1-t)r)^2.$}
\vskip 0.2truecm
\par
\noindent
By Proposition 4 $\overrightarrow{f(E)F(C)}=\overrightarrow{f(D)f(F)}$,
so $\overrightarrow{f(A)f(C)}=\overrightarrow{f(A)f(E)}+
\overrightarrow{f(E)f(C)}=\overrightarrow{f(A)f(E)}+\overrightarrow{f(D)f(F)}=
(1-t)\overrightarrow{f(A)f(D)}+(1-t)\overrightarrow{f(D)f(B)}=
(1-t)\overrightarrow{f(A)f(B)}$ and the proof is completed.
\vskip 0.2truecm
\par
It is easy to show that the present form of Theorem 3 implies a
more general form without the assumptions $0<t<1$ and $A \neq B$.
\vskip 0.2truecm
\par
\noindent
{\bf Theorem 4.} If $A,B,C,D \in {\R}^2$,
$\overrightarrow{AB}=\overrightarrow{CD}$,
$|AC|=|BD| \in \Q$ and $f:{\R}^2 \to {\F}^2$ preserves unit
distance, then $\overrightarrow{f(A)f(B)}=\overrightarrow{f(C)f(D)}$.
\vskip 0.2truecm
\par
\noindent
{\it Proof.} There exist $m \in \{0,1,2,...\}$
and $A_0,C_0,A_1,C_1,...,A_m,C_m \in {\R}^2$
such that $A_0=A$, $C_0=C$, $A_m=B$, $C_m=D$ and
for each $i \in \{0,1,...,m-1\}$ $A_iC_iC_{i+1}A_{i+1}$
is a rhombus with a rational side, see Figure 2 where $m=3$.
\vskip 0.2truecm
\par
\centerline{
\beginpicture
\setcoordinatesystem units <1mm, 1mm>
\setplotsymbol({.})
\setdashes
\plot -10 0 0 17.34 /
\plot 0 17.32 14.14 31.46 /
\plot 14.14 31.46 31.46 41.46 /
\plot 10 0 20 17.32 /
\plot 20 17.32 34.14 31.46 /
\plot 34.14 31.46 51.46 41.46 /
\plot -10 0 10 0 /
\plot 0 17.32 20 17.32 /
\plot 14.14 31.46 34.14 31.46 /
\plot 31.46 41.46 51.46 41.46 /
\setsolid
\arrow <4mm> [0.1,0.3] from -10 0 to 31.46 41.46
\arrow <4mm> [0.1,0.3] from 10 0 to 51.46 41.46
\put {$A_0=A$} at -10 -3
\put {$C_0=C$} at 10 -3
\put {$A_1$} at -3 17.32
\put {$C_1$} at 18 20
\put {$A_2$} at 9.14 31.46
\put {$C_2$} at 33 34
\put {$A_3=B$} at 31.46 44.46
\put {$C_3=D$} at 51.46 44.46
\endpicture}
\vskip 0.2truecm
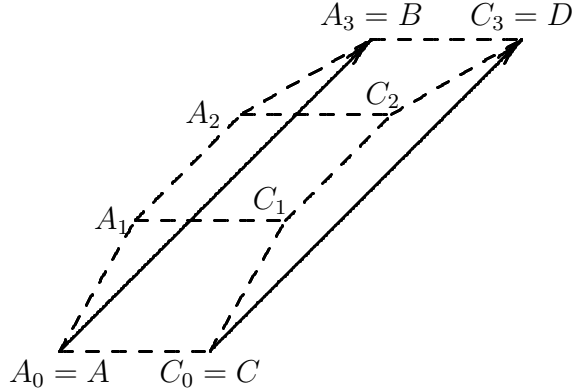
\centerline{Figure 2}
\vskip 0.2truecm
\par
\noindent
Then for each $i \in \{0,1,...,m-1\}$ $f$ preserves
the lengths of the sides of the rhombus $A_iC_iC_{i+1}A_{i+1}$.
For each $i \in \{0,1,...,m-1\}$ we have:
$\varphi_2(f(A_i),f(C_{i+1})) \neq 0$ (by Theorem 2)
and $f(C_i) \neq f(A_{i+1})$ (by Theorem 1). Therefore,
by Proposition 4 for each $i \in \{0,1,...,m-1\}$
$\overrightarrow{f(A_i)f(A_{i+1})}=\overrightarrow{f(C_i)f(C_{i+1})}$.
Hence $\overrightarrow{f(A)f(B)}=\overrightarrow{f(A_0)f(A_m)}=
\overrightarrow{f(A_0)f(A_1)}+\overrightarrow{f(A_1)f(A_2)}+...+
\overrightarrow{f(A_{m-1})f(A_m)}=\overrightarrow{f(C_0)f(C_1)}+
\overrightarrow{f(C_1)f(C_2)}+...+\overrightarrow{f(C_{m-1})f(C_m)}=
\overrightarrow{f(C)f(D)}$.
\vskip 0.2truecm
\par
\noindent
{\bf Theorem 5.} If $A,B,C,D \in {\R}^2$,
$\overrightarrow{AB}=\overrightarrow{CD}$
and $f:{\R}^2 \to {\F}^2$ preserves unit distance,
then $\overrightarrow{f(A)f(B)}=\overrightarrow{f(C)f(D)}$.
\vskip 0.2truecm
\par
\noindent
{\it Proof.} There exist $E,F \in {\R}^2$ such that
$\overrightarrow{AB}=\overrightarrow{EF}=\overrightarrow{CD}$,
$|AE|=|BF| \in {\Q}$ and $|EC|=|FD| \in {\Q}$, see Figure 3.
\vskip 0.2truecm
\par
\centerline{
\beginpicture
\setcoordinatesystem units <0.7mm, 0.7mm>
\setplotsymbol({.})
\setdashes
\plot 20 0 10 10 /
\plot 10 10 -20 0 /
\plot 30 27 20 37 /
\plot 20 37 -10 27 /
\setsolid
\arrow <4mm> [0.1,0.3] from -20 0 to -10 27
\arrow <4mm> [0.1,0.3] from 20 0 to 30 27
\arrow <4mm> [0.1,0.3] from 10 10 to 20 37
\put {$A$} at -20 -3
\put {$C$} at 20 -3
\put {$E$} at 9 6
\put {$B$} at -13 27
\put {$D$} at 33 27
\put {$F$} at 20 40
\endpicture}
\par
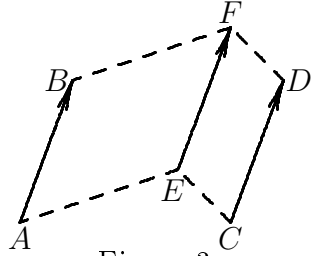
\centerline{Figure 3}
\vskip 0.2truecm
\par
\noindent
By Theorem 4
$\overrightarrow{f(A)f(B)}=\overrightarrow{f(E)f(F)}$
and
$\overrightarrow{f(E)f(F)}=\overrightarrow{f(C)f(D)}$,
so $\overrightarrow{f(A)f(B)}=\overrightarrow{f(C)f(D)}$.
\vskip 0.2truecm
\par
\noindent
As a corollary of Theorems 3 and 5 we get:
\vskip 0.2truecm
\par
\noindent
{\bf Theorem 6.} If $A,B,C,D \in {\R}^2$, $r \in \Q$,
$\overrightarrow{CD}=r\overrightarrow{AB}$ and $f:{\R}^2 \to {\F}^2$
preserves unit distance,
then $\overrightarrow{f(C)f(D)}=r\overrightarrow{f(A)f(B)}$.
\vskip 0.2truecm
\par
\noindent
{\bf Theorem 7.} If $P,Q,X,Y \in {\R}^2$,
$\overrightarrow{PQ}$ is perpendicular to $\overrightarrow{XY}$,
$|XY| \in \Q$ and $f:{\R}^2 \to {\F}^2$ preserves unit distance,
then $\overrightarrow{f(P)f(Q)}$ is perpendicular to
$\overrightarrow{f(X)f(Y)}$.
\vskip 0.2truecm
\par
\noindent
{\it Proof.} There exist $A,B,C,D,E,F \in {\R}^2$ and $r,s \in \Q$
such~that~~$\overrightarrow{PQ}=r\overrightarrow{DE}$,
$\overrightarrow{XY}=s\overrightarrow{AB}$ and the points $A,B,C,D,E,F$
form the configuration from Figure 4; it is a part of Kempe's
linkage for drawing straight lines, see \cite{Rademacher&Toeplitz}.
\vskip 0.2truecm
\par
\centerline{
\beginpicture
\setcoordinatesystem units <0.6mm,0.6mm>
\setplotsymbol({.})
\arrow <4mm> [0.1,0.3] from -40 0 to 40 0
\arrow <4mm> [0.1,0.3] from 5 18 to 5 13.2287
\plot -40 0 5 66.1437 /
\plot 35 39.6832 5 66.1437 /
\plot 40 0 35 39.6832 /
\plot 35 39.6832 5 13.2287 /
\plot 35 39.6832 5 13.2287 /
\plot 5 13.2287 20 0 /
\setdashes
\plot 5 66.1437 5 13.2287 /
\put {$A$} at -40 -4
\put {$B$} at 40 -4
\put {$F$} at 20 -4
\put {$D$} at 3 68.6437
\put {$C$} at 39 39.6862
\put {$E$} at 1.13 13.2287
\endpicture}
\vskip 0.2truecm
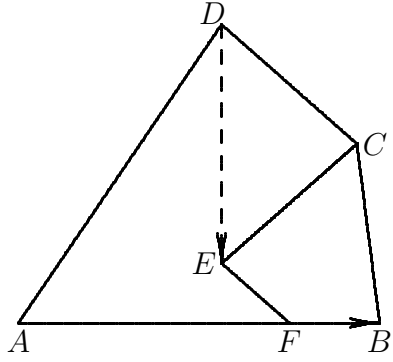
\centerline{Figure 4}
\centerline{$B \neq D$,~~$B \neq E$}
\centerline{$|AB|=|AD|=4$,~~$|CB|=|CD|=|CE|=2$,~~$|AF|=3$,~~$|FB|=|FE|=1$}
\vskip 0.4truecm
\par
\noindent
By Theorem 6 $\overrightarrow{f(P)f(Q)}=r\overrightarrow{f(D)f(E)}$
and $\overrightarrow{f(X)f(Y)}=s\overrightarrow{f(A)f(B)}$, so it
suffices to prove that
$\overrightarrow{f(D)f(E)} \cdot \overrightarrow{f(A)f(B)}=0$.
Let $a=\varphi_2(f(B),f(D))$, $b=\varphi_2(f(A),f(C))$,
$c=\varphi_2(f(B),f(E))$, $d=\varphi_2(f(C),f(F))$,
$e=\varphi_2(f(A),f(E))$.
\vskip 0.2truecm
\par
\noindent
Computing the value of $c$ we obtain:
\par
\noindent
$c=\varphi_2(f(B),f(E))=(\overrightarrow{f(B)f(E)})^2=
(\overrightarrow{f(A)f(E)}-\overrightarrow{f(A)f(B)})^2=$
\par
\noindent
$e-2\overrightarrow{f(A)f(E)} \cdot \overrightarrow{f(A)f(B)}+16$, so
\\
\centerline{~~~~~~~~~~~~~~~~~~~~~~~~~~$\overrightarrow{f(A)f(E)} \cdot \overrightarrow{f(A)f(B)}=8+\frac{1}{2}(e-c)~~~~~~~~~~~~~~~~~~~~~~~~~~~~~~~~$(4).}
\par
\noindent
Computing the value of $a$ we obtain:
\par
\noindent
$a=\varphi_2(f(B),f(D))=(\overrightarrow{f(B)f(D)})^2
=(\overrightarrow{f(A)f(D)}-\overrightarrow{f(A)f(B)})^2=$
\par
\noindent
$16-2\overrightarrow{f(A)f(D)} \cdot \overrightarrow{f(A)f(B)}+16$, so
\par
\noindent
\centerline{~~~~~~~~~~~~~~~~~~~~~~~~~~$\overrightarrow{f(A)f(D)} \cdot \overrightarrow{f(A)f(B)}=
16-\frac{1}{2}a$~~~~~~~~~~~~~~~~~~~~~~~~~~~~~~~~~~~~~(5).}
\vskip 0.2truecm
\par
\noindent
The next calculations are based
on Proposition 2 and the observation that distances $1$, $2$, $3$, $4$
are preserved by $f$.
\vskip 0.2truecm
\par
\noindent
$0=\Delta(f(A),f(B),f(E),f(F))=
\det \left[
\begin{array}{ccccc}
0 &  1 &  1 & 1 & 1 \\
1 &  0 & 16 & e & 9 \\
1 & 16 &  0 & c & 1 \\
1 &  e &  c & 0 & 1 \\
1 &  9 &  1 & 1 & 0 \\
\end{array}
\right]=
-2(e-16+3c)^2$, so
\vskip 0.2truecm
\par
\noindent
\centerline{~~~~~~~~~~~~~~~~~~~~~~~~~~~~~~~~~~~~~~~~~~~$e=16-3c$~~~~~~~~~~~~~~~~~~~~~~~~~~~~~~~~~~~~~~~~~~~~~~~(6).}
\vskip 0.2truecm
\par
\noindent
$0=\Delta(f(A),f(B),f(C),f(F))=
\det \left[
\begin{array}{ccccc}
0 &  1 &  1 &  1 &  1 \\
1 &  0 & 16 &  b &  9 \\
1 & 16 &  0 &  4 &  1 \\
1 &  b &  4 &  0 &  d \\
1 &  9 &  1 &  d &  0 \\
\end{array}
\right]
=-2(b-4d)^2,
$
so $b=4d$.
\vskip 0.4truecm
\par
\noindent
Thus $0=\Delta(f(A),f(B),f(C),f(D))=
\det \left[
\begin{array}{ccccc}
0 &  1 &  1 &  1 &  1 \\
1 &  0 & 16 &  b & 16 \\
1 & 16 &  0 &  4 &  a \\
1 &  b &  4 &  0 &  4 \\
1 & 16 &  a &  4 &  0 \\
\end{array}
\right]=$
\vskip 0.2truecm
\par
\noindent
$\det \left[
\begin{array}{ccccc}
0 &  1 &  1 &  1 &  1 \\
1 &  0 & 16 & 4d & 16 \\
1 & 16 &  0 &  4 &  a \\
1 & 4d &  4 &  0 &  4 \\
1 & 16 &  a &  4 &  0 \\
\end{array}
\right]$
$=-8a(ad+4(d^2-10d+9)).$ By Theorem 2 $a \neq 0$,
so $a=-4\frac{d^2-10d+9}{d}$.
\vskip 0.2truecm
\par
\noindent
$0=\Delta(f(B),f(C),f(E),f(F))=
\det \left[
\begin{array}{ccccc}
0 & 1 & 1 & 1 & 1 \\
1 & 0 & 4 & c & 1 \\
1 & 4 & 0 & 4 & d \\
1 & c & 4 & 0 & 1 \\
1 & 1 & d & 1 & 0 \\
\end{array}
\right]=-2c(cd+d^2-10d+9).$ By Theorem 2 $c \neq 0$, so
$c=-\frac{d^2-10d+9}{d}$. Therefore
\par
\noindent
\centerline{~~~~~~~~~~~~~~~~~~~~~~~~~~~~~~~~~~~~~~~~~~~~~~~~~$a=4c$~~~~~~~~~~~~~~~~~~~~~~~~~~~~~~~~~~~~~~~~~~~~~~~~~~~(7).}
\vskip 0.2truecm
\par
\noindent
By (4)-(7):
\vskip 0.2truecm
\par
\noindent
$\overrightarrow{f(D)f(E)} \cdot \overrightarrow{f(A)f(B)}=
(\overrightarrow{f(A)f(E)} - \overrightarrow{f(A)f(D)})
\cdot \overrightarrow{f(A)f(B)}=$
\par
\noindent
$\overrightarrow{f(A)f(E)}\cdot \overrightarrow{f(A)f(B)}-
\overrightarrow{f(A)f(D}) \cdot \overrightarrow{f(A)f(B)}=
8+\frac{1}{2}(e-c)-(16-\frac{1}{2}a)=$
\par
\noindent
$\frac{1}{2}a-\frac{1}{2}c+\frac{1}{2}e-8=
\frac{1}{2}a-\frac{1}{2}c+\frac{1}{2}(16-3c)-8=
\frac{1}{2}a-2c=0$.
\vskip 0.2truecm
\par
\noindent
The proof is completed.
\vskip 0.2truecm
\par
\noindent
{\bf Theorem 8.} If $A,B,C,D \in {\R}^2$, $\overrightarrow{AB}$
and $\overrightarrow{CD}$ are linearly dependent and
$f:{\R}^2 \to {\F}^2$ preserves unit distance,
then $\overrightarrow{f(A)f(B)}$ and 
$\overrightarrow{f(C)f(D)}$ are linearly dependent.
\vskip 0.2truecm
\par
\noindent
{\it Proof.} We choose $X,Y \in {\R}^2$ such that $|XY|=1$ and
both vectors $\overrightarrow{AB}$ and $\overrightarrow{CD}$
are perpendicular to $\overrightarrow{XY}$.
Then $\varphi_2(f(X),f(Y))=1$, so obviously $f(X) \neq f(Y)$.
By Theorem 7 both vectors $\overrightarrow{f(A)f(B)}$
and $\overrightarrow{f(C)f(D)}$ are perpendicular
to $\overrightarrow{f(X)f(Y)}$. These two facts
imply that $\overrightarrow{f(A)f(B)}$ and
$\overrightarrow{f(C)f(D)}$ are linearly dependent.
\vskip 0.2truecm
\par
\noindent
As a corollary of Theorem 8 we get:
\vskip 0.2truecm
\par
\noindent
{\bf Theorem 9.} Unit-distance preserving mappings
from ${\R}^2$ to ${\F}^2$ preserve collinearity of points.
\vskip 0.2truecm
\par
Now, we describe a general form of a unit-distance preserving
mapping from ${\R}^2$ to~${\F}^2$.
Assume that $f:{\R}^2 \to {\F}^2$ preserves unit distance.
By Theorem 5 the mapping
\vskip 0.2truecm
\par
\noindent
\centerline{${\R}^2 \ni u \stackrel{\varphi}{\longrightarrow} f(u)-f(0) \in {\F}^2$}
\par
\noindent
satisfies
\vskip 0.2truecm
\par
\noindent
~~~~~~~~~~~~~~~~~~~~~~~~~~~~~~~~$\forall u,v \in {\R}^2~~ \varphi(u+v)=\varphi(u)+\varphi(v)$~~~~~~~~~~~~~~~~~~~~~~~~~~~~~~~~~(8).
\vskip 0.2truecm
\par
\noindent
By Theorem 1 for each non-zero $u \in {\R}^2$
$\varphi(u) \neq 0$, obviously $\varphi(0)=0$.
By Theorem~9 $\varphi$ preserves collinearity of points.
Therefore, for each $\lambda \in \R$ and each non-zero
$u \in {\R}^2$ there exists unique $\theta(\lambda,u) \in \F$ such
that $\varphi(\lambda u)=\theta(\lambda,u)\varphi(u)$.
It follows from this and~(8) that $\varphi$ transforms linearly
independent vectors into linearly independent vectors.
Otherwise it would imply that $f({\R}^2)$ is contained in an affine
line. It is impossible because $f({\R}^2)$ contains points $c_1$, $c_2$, $c_3$
such that $\varphi_2(c_1,c_2)=\varphi_2(c_1,c_3)=\varphi_2(c_2,c_3)=1$, so
$$
\Delta(c_1,c_2,c_3)=
\det \left[
\begin{array}{cccc}
0 & 1 & 1 & 1 \\
1 & 0 & 1 & 1 \\
1 & 1 & 0 & 1 \\
1 & 1 & 1 & 0 \\
\end{array}\;\right]=-3 \neq 0$$
and the points $c_1$, $c_2$, $c_3$ are affinely independent by Proposition 1.
\vskip 0.2truecm
\par
Now we prove that $\theta(\lambda,u)$ does not depend on $u$,
the exposition follows the one in~\cite{Lelong-Ferrand}.
We consider two cases:
\vskip 0.2truecm
\par
\noindent
{\bf (I)} Non-zero vectors $u,v \in {\R}^2$ are linearly
independent. Then for each $\lambda \in \R$ we have:
\par
\noindent
\centerline{$\theta(\lambda,u+v)\varphi(u)+\theta(\lambda,u+v)\varphi(v)=
\theta(\lambda,u+v)\varphi(u+v)=$}
\centerline{$\varphi(\lambda(u+v))=
\varphi(\lambda u)+\varphi(\lambda v)=
\theta(\lambda,u)\varphi(u)+\theta(\lambda,v)\varphi(v).$}
\vskip 0.2truecm
\par
\noindent
Since $\varphi(u)$ and $\varphi(v)$ are linearly indepenent,
we conclude that 
$\theta(\lambda,u+v)=
\theta(\lambda,u)$
and $\theta(\lambda,u+v)=\theta(\lambda,v)$. Hence
$\theta(\lambda,u)=\theta(\lambda,v)$.
\vskip 0.2truecm
\par
\noindent
{\bf (II)} Non-zero vectors $u,v \in {\R}^2$ are linearly
dependent. Then there exists $w \in {\R}^2$ such that $u$, $w$
are linearly independent and $v,w$ are linearly independent.
Applying {\bf (I)} we get:
\par
\centerline{$\forall \lambda \in \R$
~
$\theta(\lambda,u)=\theta(\lambda,w)=\theta(\lambda,v)$.}
\vskip 0.2truecm
\par
Let $\rho: \R \to \F$, $\rho(\lambda)=\theta(\lambda,u)$;
we have proved that $\theta(\lambda,u)$ does not depend on $u$.
Now we prove that $\rho$ is a field homomorphism,
the exposition follows the one in \cite{Lelong-Ferrand}.
For each $\lambda, \mu \in \R$ and each non-zero $u \in {\R}^2$
we have:
$$
\rho(\lambda+\mu)\varphi(u)=
\varphi((\lambda+\mu)u)=
\varphi(\lambda u)+\varphi(\mu u)=
\rho(\lambda)\varphi(u)+\rho(\mu)\varphi(u)=
(\rho(\lambda)+\rho(\mu))\varphi(u).
$$
\par
\noindent
Therefore
\par
\noindent
\vskip 0.2truecm
\centerline{$\forall \lambda, \mu \in \R$
~
$\rho(\lambda+\mu)=\rho(\lambda)+\rho(\mu)$.}
\vskip 0.2truecm
\par
\noindent
For each $\lambda, \mu \in \R$ and each non-zero $u \in {\R}^2$
we have:
$$
\rho(\lambda \mu)\varphi(u)=
\varphi((\lambda \mu)u)=
\varphi(\lambda(\mu u))=
\rho(\lambda)\varphi(\mu u)=
\rho(\lambda)\rho(\mu)\varphi(u).
$$
\par
\noindent
Therefore
\vskip 0.2truecm
\par
\noindent
\centerline{
$\forall \lambda, \mu \in \R$
~
$\rho(\lambda \mu)=\rho(\lambda)\rho(\mu)$.}
\vskip 0.2truecm
\par
\noindent
We have proved that $\rho$ is a field homomorphism.
It gives our main result:
\vskip 0.2truecm
\par
\noindent
{\bf Theorem 10.} Each unit-distance preserving mapping
$f:{\R}^2 \to {\F}^2$ has a form $I \circ (\rho,\rho)$,
where $\rho: \R \to \F$ is a field homomorphism and
$I:{\F}^2 \to {\F}^2$ is an affine mapping with orthogonal
linear part.
\vskip 0.2truecm
\par
{\bf Acknowledgement.} The author wishes to thank the anonymous reviewer
for providing necessary assumptions on the field $\F$. The newest version
of \cite{Tyszka200??} is available on
{\it http://arxiv.org/abs/math.MG/0302276}.

Apoloniusz Tyszka\\
Technical Faculty\\
Hugo Ko\l{}\l{}\k{a}taj University\\
Balicka 104, 30-149 Krak\'ow, Poland\\
E-mail address: {\it rttyszka@cyf-kr.edu.pl}
\end{document}